\documentclass[conference]{IEEEtran}
\usepackage{times}

\usepackage[left=15.9mm, right=15.9mm,top=19.1mm,bottom=26mm]{geometry}

\usepackage[numbers]{natbib}
\usepackage{multicol}
\usepackage[bookmarks=true]{hyperref}
\usepackage[cmex10]{amsmath}

\usepackage{graphicx}
\usepackage{acronym}
\usepackage{xspace}
\usepackage{pgfplots}
\usepackage[noend]{algpseudocode}
\usepackage{algorithm}
\usepackage{algorithmicx}
\usepackage{xcolor}
\usepackage{cancel}

\algrenewcommand\algorithmicrequire{\textbf{Input:}}
\algrenewcommand\algorithmicensure{\textbf{Input:}}

\DeclareMathOperator{\prox}{prox}

\DeclareMathOperator*{\argmin}{argmin}

\definecolor{amethyst}{rgb}{1, 0, 1}
\definecolor{blue-violet}{rgb}{0.54, 0.17, 0.89}
\definecolor{brightturquoise}{rgb}{0.03, 0.91, 0.87}
\usepackage{pgfplots}\pgfplotsset{compat=1.4}

\usepackage{caption}
\usepackage{subcaption}

\newcommand{\blue}{\textcolor{black}}

\begin{document}

\title{Robust Singular Smoothers For Tracking Using Low-Fidelity Data}




%
\author{\authorblockN{Jonathan Jonker\authorrefmark{1},
Aleksandr Aravkin\authorrefmark{2},
James V. Burke\authorrefmark{1}, 
Gianluigi Pillonetto\authorrefmark{3} and
Sarah Webster\authorrefmark{4}}
\authorblockA{\authorrefmark{1}Department of Mathematics, University of Washington, Seattle, WA USA, \it{jonkerjo@uw.edu, jvburke@uw.edu}}
\authorblockA{\authorrefmark{2}Department of Applied Mathematics, University of Washington, Seattle, WA USA, 
\it{saravkin@uw.edu} }
\authorblockA{\authorrefmark{3}Department of Information Engineering, University of Padova, Italy, 
\it{giapi@dei.unipd.it} }
\authorblockA{\authorrefmark{4}Applied Physics Laboratory, University of Washington, Seattle, WA USA, 
{\it swebster@apl.washington.edu} } }

\maketitle

\begin{abstract}
Tracking underwater autonomous platforms is often difficult because of noisy, biased, and discretized input data. Classic filters and smoothers based on standard assumptions of Gaussian white noise break down when presented with any of these challenges.
Robust models (such as the Huber loss) and constraints (e.g. maximum velocity) are used to attenuate these issues. 
Here, we consider robust smoothing with singular covariance, which covers bias and correlated noise, as well as many specific model types, such as those used in navigation.  
In particular, we show how to combine singular covariance models with robust losses and state-space constraints in a unified framework that can handle very low-fidelity data. A noisy, biased, and discretized navigation dataset from a submerged, low-cost inertial measurement unit (IMU) package, with ultra short baseline (USBL) data for ground truth, provides an opportunity to stress-test the proposed framework with promising results.  We show how robust modeling elements improve our ability to analyze the data, and  present batch processing results 
for 10 minutes of data with three different frequencies of available USBL position fixes
(gaps of 30 seconds, 1 minute, and 2 minutes). 
The results suggest that the framework can be extended to real-time tracking 
using robust windowed estimation. 

\end{abstract}

\IEEEpeerreviewmaketitle

\section{Introduction}


State-space models are ubiquitous in signal processing, and allow integration of disparate measurements to inform estimation, decisions, and control. Classic filtering~\cite{kalman1960new,KalBuc} and smoothing~\cite{rauch1965maximum,fraser1969optimum} are core tools used to estimate these models. Their dependence on high-fidelity data, driven by Gaussian assumptions on errors and innovations,
has demanded unequivocal attention from researchers and practitioners, and inspired robust dynamic inference methods. 

While early robust approaches \citep{Kassam1985,Schick1994} sought to modify iterations of the Kalman filter (KF) and Rauch-Tung-Striebel (RTS) smoother, over the last 25 years researchers have used {\it robust formulations} to weave 
assumptions on errors and innovations directly into the estimation problems themselves~\cite{Fahr1991,Bell1994,kim2009ell_1,aravkin2011ell,Farahmand2011,aravkin2013sparse,aravkin2014robust}. Constraints, when available, are also readily incorporated into the problem formulation~\cite{bell2009inequality}.  
Specifying the formulation leaves one free to choose from a range of optimization algorithms; the survey~\cite{aravkin2017generalized} describes a general class of models as well as first- and second-order methods to solve them. 

Our focus is on models with singular variances for process and measurement residuals. These  models are excluded by the assumptions of generalized smoothing~\cite{aravkin2017generalized} and all of the various special cases cited in the survey. 
In this paper, we build on the recently proposed framework of~\cite{jonker2019fast} for singular models,
and systematically develop complementary modeling elements: {\bf robust penalties, informative constraints, and singular models}. 
 The resulting approach exploits the structure of singular covariances  head on rather than using workarounds such as pseudo-inverses or variance boosting that either 
do not work in the general setting or introduce unnecessary changes to the fundamental model (see discussion in~\cite{jonker2019fast}).
\begin{figure}[h!]
\begin{center}
\includegraphics[scale=.5]{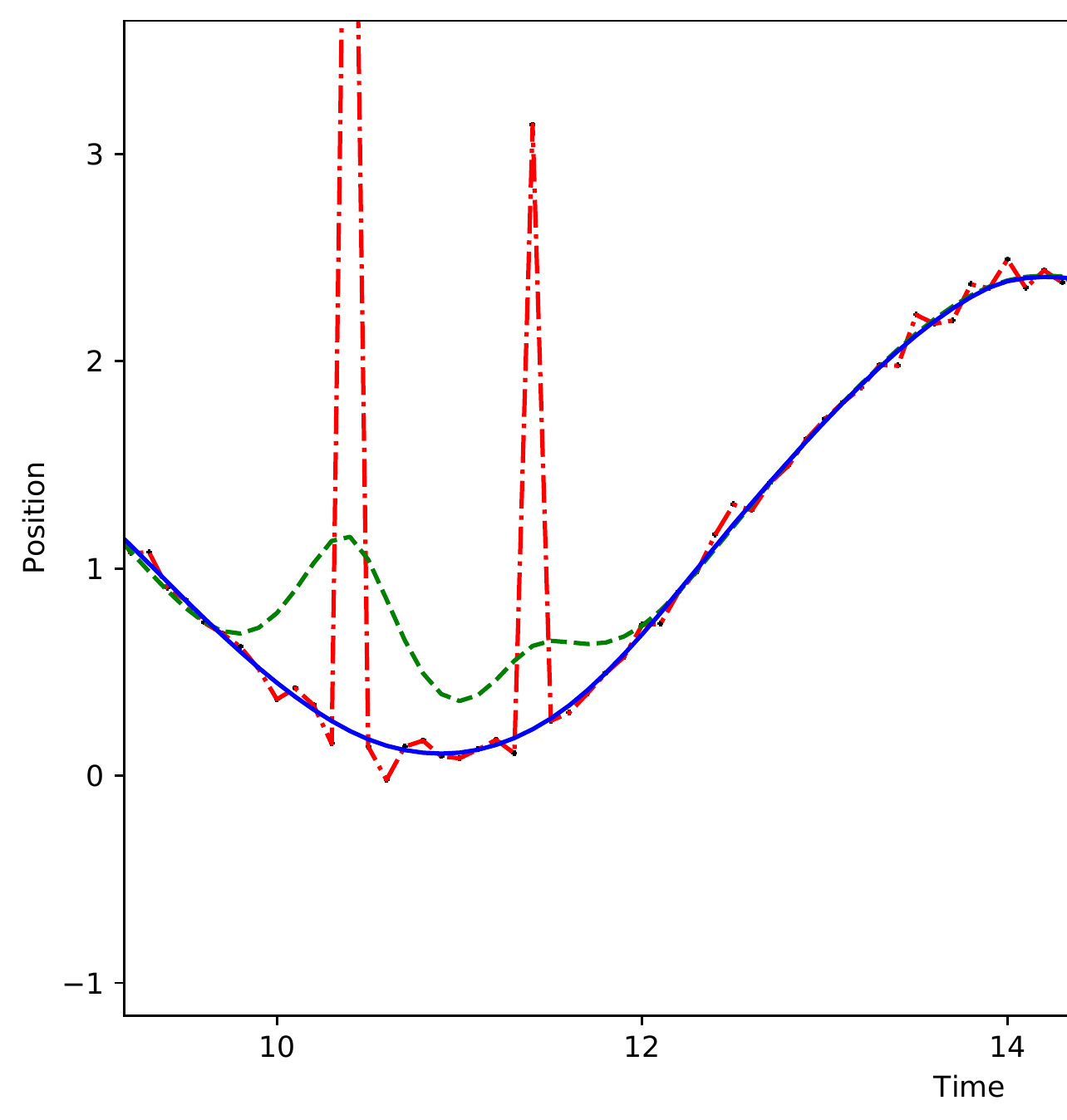}
\caption{\label{fig:intro} We track a simple trajectory in the presence of outliers. The red dash-dot shows a`robust' Huberized approach implemented using a pseudo-inverse; green dash  shows the proposed singular $\ell^2$ estimate; blue solid shows the proposed singular Huber estimate, which clearly tracks the true state. }
\end{center}
\end{figure}
A simple synthetic tracking example with a singular process model shows how the common tack of replacing the inverse by the pseudo-inverse fails dramatically in the presence of outliers (Figure~\ref{fig:intro}). Any robust approach requires control of the null spaces associated to process and observations.  We develop a direct practical formulation and method, and test it on a batch smoothing analysis of real field data, with a view towards real-time implementation in future work. 


\noindent
{\bf Background.} 
Our main goal is to infer an unobserved state sequence $x_1, \dots, x_N$ from 
noisy observations $y_1, \dots, y_N$ using the model: 
\begin{equation} \label{eq:statespace}
\begin{aligned}
x_1 & =  x_0 + w_1
\\
x_k & =  G_kx_{k-1} + w_k, \quad  k = 2, \dots, N
\\
y_k & =  H_kx_k + v_k, \quad  k = 1, \dots, N,\\
x_k &\in X_k, \quad \mbox{each $X_k$ polyhedral,}
\end{aligned}
\end{equation}
where $x_0$ is a given initial state estimate, $G_k$ and $H_k$ are linear process and measurement models, 
$y_1, \dots, y_N$ are observations, and $X_k$ specify additional information through constraints. 
The framework of~\cite{aravkin2017generalized}
assumes that $w_k$ and $v_k$ are mutually independent random variables with known {\it nonsingular} covariances $Q_k$ and $R_k$, and that they follow from log-concave distributions; in particular they may be non-Gaussian. 

Synthesizing all of this information gives the problem
\begin{equation}
    \label{eq:original}
    \begin{aligned}
    \min_{x_1\in X_1, \dots, x_n \in X_N} &\sum_k\rho_p(Q_k^{-1/2}(G_kx_k-x_{k-1})) \\
     &+ \rho_m(R_k^{-1/2}(H_kx_k- y_k)),
    \end{aligned}
\end{equation}
with $\rho_p$ and $\rho_m$ convex penalties. Using~\eqref{eq:original} provides estimates that 
are robust to outliers and can follow sudden changes in the state. Most of the inference- 
or optimization-based work in the convex dynamic setting is a special case of~\eqref{eq:original}.
Many examples, including robust penalties and constraints, are collected in~\cite{aravkin2017generalized}.

We now extend to {\it singular covariances} $R_k$ (for errors $v_k$) and $Q_k$ 
(for innovations $w_k$). These models specify key use cases, particularly for innovations (process) modeling, 
briefly summaried below (see~\cite{jonker2019fast} for a more  detailed discussion).

\noindent
{\bf Deterministic integrals.} Most models in robotics, particularly in navigation, use integration to model process relationships between state variables (e.g. when position, velocity, and acceleration are part of the state). Any deterministic integral yields a singular process model. The simplest example 
(with position a direct integral of stochastic velocity) is used to create Figure~\ref{fig:intro}.

\noindent
{\bf Nuisance parameters.} Unknown constants that need calibration (such as fixed 
instrument biases) require special modeling in the nonsingular 
paradigm~\eqref{eq:original}. With singular models, we can augment the state in order to infer these parameters. 

\noindent
{\bf Auto-regressive models and correlated errors.} 
State-space models are broadly used in auto-regressive, moving average, 
and time series models~\cite{hyndman2008forecasting}. These elements 
also appear in general smoothing models, 
particularly to deal with correlated measurement errors~\cite{chui2017kalman}.

All three examples are accessible in the classic linear Gaussian setting. The KF need not invert $Q$ or $R$, 
and provides the minimum variance estimate for both singular and nonsingular models~\cite{ansley1982geometrical}. 
Some algorithms rely on precise knowledge of the error structure or explicit equality constraints~\cite{koopman1997exact,ko2007state,ait2011fixed}. 
Correlated errors are dealt with by augmenting the state 
and using a singular model~\cite{chui2017kalman}. 
{\bf None of these techniques generalize to singular models in the setting of~\eqref{eq:original}}, and some naive generalizations fail dramatically (Figure~\ref{fig:intro}).  

This paper builds on the reformulation of~\cite{jonker2019fast} for singular models. 
We develop a systematic approach and test it using a navigation model with a real-world mooring dataset. We show how robust statistics, singular models, and constraints can be systematically used to overcome a range of challenges simultaneously present in the dataset: (1) outliers, (2) deterministic relationships between states, (3) measurement biases, and (4) coarsely discretized observations.  

The paper proceeds as follows. Section~\ref{sec:summary} summarizes 
the singular formulation of~\cite{jonker2019fast} and relevant optimization algorithms. Section~\ref{sec:modeling} develops the key modeling elements to address common data challenges.  Section~\ref{sec:navmodels} presents the navigation models. 
Section~\ref{sec:analysis} shows how the model elements come together to analyze the target mooring dataset, and obtain a high 
fidelity track from low-fidelity observations. 
Section~\ref{sec:conclusion} concludes with discussion and future work.

\section{Robust Singular Formulation and Algorithm}
\label{sec:summary}

We reformulate the robust smoothing problem to 
seamlessly allow both singular and nonsingular covarinace models for errors and innovations. 
The resulting problem can be solved with any primal-dual algorithm. 
We show how to apply the classic Douglas-Rachford splitting (DRS) algorithm 
(see e.g.~\cite{eckstein1992douglas,davis2016convergence}) to the reformulated problem.

Problem~\eqref{eq:original} can be reformulated by introducing auxiliary variables $u_k$ 
and $t_k$ to represent pre-whitened innovations and measurement residuals: 
\begin{equation}
    \label{eq:reform}
        \begin{aligned}
    \min_{x, u, t} & \quad \sum_k\rho_p(u_k) + \rho_m(t_k) + \rho_s(x_k) \\
   \text{s.t.} & \quad  Q_k^{1/2}u_k = G_kx_k - x_{k-1}\\
   & \quad R_k^{1/2}t_k = H_kx_k-y_k
\end{aligned}
\end{equation}
where $\rho_s(x_k)$ may be taken as the convex indicator function to recover the constraints 
in~\eqref{eq:original}:
\[
\rho_s(x_k) = 
\begin{cases}
0 \quad  x_k \in X_k\\
\infty \quad x_k \not \in X_k.
\end{cases}
\]
When $Q_k$ and $R_k$ are invertible, we can solve for $u_k,t_k$ and recover~(\ref{eq:original}). 
Otherwise, problem~(\ref{eq:reform}) is well-posed while~(\ref{eq:original}) is not. 
We can write~\eqref{eq:reform} in compact form 
\begin{equation}
\label{eq:full}
\begin{aligned}
\min_{z}  &\quad \rho(z)  \quad \text{s.t. } Az = \hat{w}, \\
 \rho(z)  &= \sum_{k=1}^N \rho_p(u_k) + \rho_m(t_k) + \rho_s(x_k).
\end{aligned}
\end{equation}
where 
\begin{equation}
\label{eq:zw}
\begin{aligned}
z^T &= \begin{pmatrix} u_1^T & t_1^T & x_1^T & \dots u_N^T & t_N^T & x_N^T\end{pmatrix} \\
\hat{w}^T &= \begin{pmatrix} x_0^T & y_1^T & 0 & y_2^T & \dots & 0 & y_N^T\end{pmatrix},
\end{aligned}
\end{equation}
\begin{equation}
\label{eq:A}
A = 
\begin{pmatrix} D_1 & 0 & \dots & 0 \\ 
B_1 & D_2 & 0 & \vdots \\ 
 0& \ddots & \ddots &  0\\
  0 & 0 & B_{N-1} & D_N\end{pmatrix},
\end{equation}
and 
\begin{equation*}
\label{eq:DjBj}
\begin{aligned}
D_i = \begin{pmatrix}Q_i^{1/2} & 0 & I\\ 0 & R_i^{1/2} & H_i\end{pmatrix}, 
B_j = \begin{pmatrix}0 & \qquad 0 & -G_{j+1}\\ 0 & \qquad 0 &  0\end{pmatrix}.
\end{aligned}
\end{equation*}

The variables are ordered in such a way that $A$ is block bi-diagonal. 
If all observations $z_i$ lie in the range of $H_i$, the constraint $Az = \hat{w}$ will 
be feasible~\cite{jonker2019fast}.

The problem~\eqref{eq:full} is a convex optimization problem and can be solved using a 
variety of techniques. We show that the DRS algorithm is straightforward to implement,
and preserves the computational complexity of the classic KF/RTS algorithms 
because of the way $A$ is structured in~\eqref{eq:A}.

Given a convex function $f$, 
its convex conjugate $f^*$ is given by
\[
f^*(y) = \sup_x \langle x, y \rangle - f(x),
\]
and its proximal operator with step $\alpha$, 
denoted by $\prox_{\alpha f}$ (see e.g. \cite{combettes2011proximal}) is given by:
\begin{equation}
    \label{eq:prox}
    \prox_{\alpha f}(\zeta) = \arg\min_x \frac{1}{2\alpha} \|\zeta - x\|^2 + f(x).
\end{equation}

For a long list of objectives, prox operators
are available in closed form  or  are efficiently computable.
In particular this is the case when $\rho_p, \rho_m, \rho_s$ form any subset of the numerous elements briefly surveyed in Section~\ref{sec:modeling}.   
It is actually the prox of $\rho^*$ that appears in the DRS iteration (Algorithm~\ref{alg:DRS}) rather than the prox of $\rho$, but these are linked by the simple formula 
\[
\prox_\rho(z) + \prox_{\rho^*}(z)= z. 
\]

To specify the algorithm, we let $g(z)$ be the indicator of the affine feasible  region $Az = \hat w$:
\[
g(z) = 
\begin{cases}
0 \quad Az=\hat{w}\\
\infty \quad Az \neq \hat{w}
\end{cases}
\]

Problem~\eqref{eq:full} can now be written simply as 
\[
\min_{z} \rho(z) + g(z)
\]
which is a natural template for DRS, detailed in Algorithm~\ref{alg:DRS}. 
    \begin{algorithm}[H]
  \caption[Caption]{Douglas-Rachford Splitting (DRS)
    \label{alg:DRS}}
  \begin{algorithmic}[1]
    \Require{Initialize at any $z^0$, $\zeta^0$.}
    \Loop
    \State {$z^k = \prox_{\tau g}(z^{k-1}-\tau \zeta^{k-1})$} 
     \State {$\zeta^k = \prox_{\sigma \rho^*}(\zeta^{k-1}+\sigma(2z^{k} - z^{k-1}))$}  
       \EndLoop
       \Return{$z^k$}
  \end{algorithmic}
\end{algorithm}

To implement Algorithm~\ref{alg:DRS} we need proximal operators of $\rho_p, \rho_m$, and $\rho_s$. Eight common piecewise linear-quadratic (PLQ) penalties are shown in Figure~\ref{fig:PLQ}, and their proximal operators
are summarized in Table~\ref{table:PLQ}.

The proximal operator for $g$ is given by
\[
\prox_g(\eta) = \argmin_{Az=\hat{w}} \frac{1}{2}||\eta - z||^2
\]
which is a least squares problem with affine constraints.
Solving it efficiently leverages the structure of~\eqref{eq:A}. In particular we need to solve 
a single structured linear system
\begin{equation}
    \label{eq:linsystem}
\begin{bmatrix}I & A^T\\ 0 & AA^T\end{bmatrix}\begin{bmatrix}z\\ \nu \end{bmatrix} = \begin{bmatrix}\eta \\ A\eta - \hat{w}\end{bmatrix}
\end{equation}
where $AA^T$ is block tridiagonal and does not change between iterations. In our implementation, we need only compute a single block bidiagonal factorization once, which can then be used to solve~(\ref{eq:linsystem}) in $O(n^2N)$ operations 
in each iteration, no more expensive than a single matrix-vector multiply.

For piecewise-linear quadratic $\rho$~\cite{rockafellar2009variational,aravkin2017generalized}, 
DRS converges to an optimal solution at a local linear rate~\cite{jonker2019fast}, which does not depend on the condition number of $A$. 
A good initialization makes DRS competitive with the fastest
available solvers, even second order methods with quadratic local rates~\cite{aravkin2017generalized}.

\section{Modeling Elements}
\label{sec:modeling}

\begin{figure}[t!]
    \begin{subfigure}[t]{0.24\textwidth}
       \centering
\begin{tikzpicture}
  \begin{axis}[
    thick,
    height=2cm,
    xmin=-2,xmax=2,ymin=0,ymax=1,
    no markers,
    samples=50,
    axis lines*=left, 
    axis lines*=middle, 
    scale only axis,
    xtick={-1,1},
    xticklabels={},
    ytick={0},
    ] 
\addplot[blue, domain=-2:+2]{.5*x^2};
  \end{axis}
  \end{tikzpicture}
  \caption{\label{fig:quadratic}quadratic}
    \end{subfigure}%
    \begin{subfigure}[t]{0.24\textwidth}
        \centering
\begin{tikzpicture}
  \begin{axis}[
    thick,
    height=2cm,
    xmin=-2,xmax=2,ymin=0,ymax=1,
    no markers,
    samples=100,
    axis lines*=left, 
    axis lines*=middle, 
    scale only axis,
    xtick={-1,1},
    xticklabels={},
    ytick={0},
    ] 
  \addplot[red, dashed, domain=-2:+2]{abs(x)};
  \end{axis}
\end{tikzpicture}
    \caption{\label{fig:1norm}1-norm}
    \end{subfigure}
    \begin{subfigure}[t]{0.24\textwidth}
        \centering
\begin{tikzpicture}
  \begin{axis}[
    thick,
    height=2cm,
    xmin=-2,xmax=2,ymin=0,ymax=1,
    no markers,
    samples=50,
    axis lines*=left, 
    axis lines*=middle, 
    scale only axis,
    xtick={-1,1},
   xticklabels={},
    ytick={0},
    ] 
\addplot[red,domain=-2:0,densely dashed]{-.3*x};
\addplot[red,domain=0:+2,densely dashed]{.7*x};
  \end{axis}
\end{tikzpicture}
            \caption{\label{fig:quantile}quantile, $\tau = 0.3$}
        \end{subfigure}
   \begin{subfigure}[t]{0.24\textwidth}
   \centering
   \begin{tikzpicture}
  \begin{axis}[
    thick,
    height=2cm,
    xmin=-2,xmax=2,ymin=0,ymax=1,
    no markers,
    samples=50,
    axis lines*=left, 
    axis lines*=middle, 
    scale only axis,
    xtick={-1,1},
    xticklabels={},
    ytick={0},
    ] 
\addplot[red,domain=-2:-1,densely dashed]{-x-.5};
\addplot[blue, domain=-1:+1]{.5*x^2};
\addplot[red,domain=+1:+2,densely dashed]{x-.5};
\addplot[blue,mark=*,only marks] coordinates {(-1,.5) (1,.5)};
  \end{axis}
\end{tikzpicture}
\caption{\label{fig:huber}huber, $\kappa = 1$}
\end{subfigure}    
   \begin{subfigure}[t]{0.24\textwidth}
   \centering
\begin{tikzpicture}
  \begin{axis}[
    thick,
    height=2cm,
    xmin=-2,xmax=2,ymin=0,ymax=1,
    no markers,
    samples=100,
    axis lines*=left, 
    axis lines*=middle, 
    scale only axis,
    xtick={-.24,.56},
    xticklabels={},
    ytick={0},
    ] 
\addplot[red,domain=-2:-2*0.3*0.4,densely dashed]{0.3*abs(x) - 0.4*0.3^2};
\addplot[blue,domain=-2*0.3*0.4:2*(1-0.3)*0.4]{0.25*x^2/0.4};
\addplot[red,domain=2*(1-0.3)*0.4:2,densely dashed]{(1-0.3)*abs(x) - 0.4*(1-0.3)^2};
\addplot[blue,mark=*,only marks] coordinates {(-.24,0.0550) (0.56,0.20)};
  \end{axis}
\end{tikzpicture}
\caption{\label{fig:huberQ} quantile huber}
\end{subfigure}  
   \begin{subfigure}[t]{0.24\textwidth}
   \centering
\begin{tikzpicture}
  \begin{axis}[
    thick,
    height=2cm,
    xmin=-2,xmax=2,ymin=0,ymax=1,
    no markers,
    samples=50,
    axis lines*=left, 
    axis lines*=middle, 
    scale only axis,
    xtick={-0.5,0.5},
    xticklabels={},
    ytick={0},
    ] 
    \addplot[red,domain=-2:-0.5,densely dashed] {-x-0.5};
    \addplot[domain=-0.5:+0.5] {0};
    \addplot[red,domain=+0.5:+2,densely dashed] {x-0.5};
    \addplot[blue,mark=*,only marks] coordinates {(-0.5,0) (0.5,0)};
  \end{axis}
\end{tikzpicture}
\caption{\label{fig:vapnik} vapnik, $\epsilon = 0.5$}
\end{subfigure}  
   \begin{subfigure}[t]{0.24\textwidth}
   \centering
\begin{tikzpicture}
  \begin{axis}[
    thick,
    height=2cm,
    xmin=-2,xmax=2,ymin=0,ymax=1,
    no markers,
    samples=50,
    axis lines*=left, 
    axis lines*=middle, 
    scale only axis,
    xtick={-1,1, -.5, .5},
    xticklabels={},
    ytick={0},
    ] 
\addplot[domain=-0.25:+0.25] {0};
\addplot[red,domain=-2:-1,densely dashed]{-x-.5-.5*.25};
\addplot[blue, domain=-1:-.25]{.5*x^2-.5*.25};
\addplot[blue, domain=.25:1]{.5*x^2-.5*.25};
\addplot[red,domain=+1:+2,densely dashed]{x-.5-.5*.25};
\addplot[blue,mark=*,only marks] coordinates {(-1,.5-.5*.25) (1,.5-.5*.25)(-.45, 0) (.45, 0)};
  \end{axis}
\end{tikzpicture}
\caption{\label{fig:sel} hubnick}
\end{subfigure}  
   \begin{subfigure}[t]{0.24\textwidth}
   \centering
\begin{tikzpicture}
  \begin{axis}[
    thick,
    height=2cm,
    xmin=-2,xmax=2,ymin=0,ymax=1,
    no markers,
    samples=100,
    axis lines*=left, 
    axis lines*=middle, 
    scale only axis,
    xtick={-1,1},
    xticklabels={},
    ytick={0},
    ] 
\addplot[amethyst, domain=-2:+2]{.5*x^2 + 0.5*abs(x)};
  \end{axis}
\end{tikzpicture}
\caption{\label{fig:enet} elastic net}
\end{subfigure}  
    \caption{\label{fig:PLQ}Common piecewise linear-quadratic (PLQ) losses.}
\end{figure}
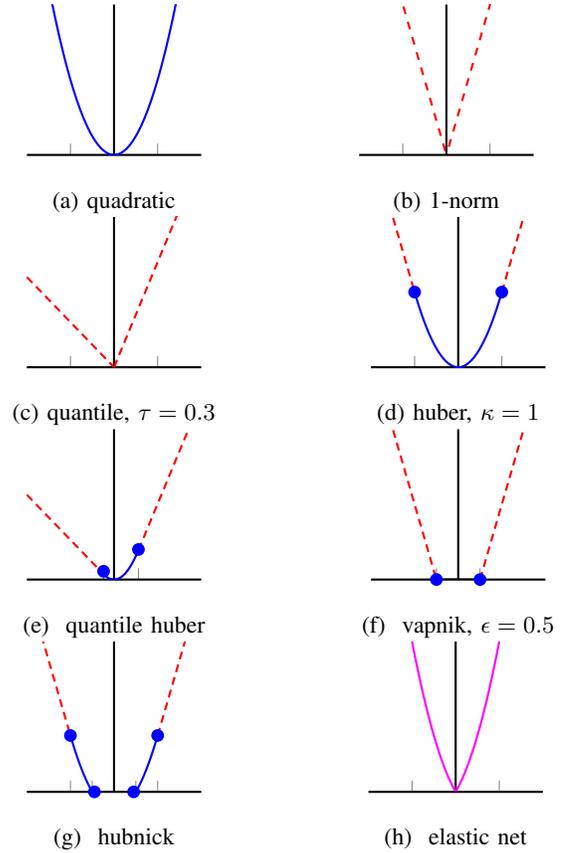

The proposed framework has three complementary modeling elements:
singular covariance matrices $Q$ and $R$;  process/measurement penalties $\rho_p, \rho_m$; and constraints 
$\rho_s$ on the state. 
In this section, we show a range of choices for each element,
and compute the operators required for Algorithm~\ref{alg:DRS}.  

\noindent
{\bf Singular covariances} 
can be used to capture affine constraints, auto-regressive structure, integrated errors, and bias. 
\begin{itemize}
\item \emph{Affine constraints using singular $R$.}  
the $i$th element of the state 
at time $k$ is known exactly, add row 
\[
\begin{bmatrix} 0 & \dots &0 & \underbrace{1}_{i} & 0 &\dots &0\end{bmatrix}
\]
to the measurement model $H_k$, a row and column of zeros to $R_k$, and the known value as the last element of $z_k$.

\item \emph{Bias with singular $Q$}. A common model for bias is to include it as a non-varying component 
across the state:
\[
\widetilde x_k = \begin{bmatrix} x_k \\ b \end{bmatrix}, \quad \widetilde Q_k = \begin{bmatrix} Q & 0\\
0 & 0 \end{bmatrix}.
\]

\item \emph{Correlated noise using singular $Q$.}
Correlated noise $w_k$ is typically modeled by~\cite{Chui2009} 
\[
w_k = M w_{k-1} + \beta_k, \quad \beta_k \sim N(0, Q). 
\]
Here too we can augment the state and use a singular process variance:
\[
\widetilde x_k = \begin{bmatrix} x_k \\ w_k \end{bmatrix}, 
\widetilde G_k = \begin{bmatrix} G_k & I \\ 0 & M \end{bmatrix}, 
\quad \widetilde Q_k = \begin{bmatrix} 0 & 0 \\
0 & Q_k \end{bmatrix}.
\]

\end{itemize}

\noindent
{\bf Piecewise linear-quadratic (PLQ) Penalties.}
The proposed framework allows process innovations, measurement 
residuals, and state regularization to come from any convex prox-friendly 
penalty. To keep the exposition simple, we collect eight commonly used
convex piecewise linear-quadratic penalties in Figure~\ref{fig:PLQ}, 
and compute their prox operators in Table~\ref{table:PLQ}.
The penalties can be thought of in terms of three features: 
\begin{itemize}
\item Behavior at origin: nonsmooth features encourage 
exact fitting of the quantity being measured, while deadzones 
are appropriate for discretized observations. 
\item Tail growth: asymptotically linear penalties are more tolerant of large inputs. Applied to measurements, this gives robustness to outliers; applied to innovations, it gives an ability to quickly track evolving trends. 
\item Asymmetry: allows handling of special cases where under-estimating is qualitatively different from over-estimating.
\end{itemize}

\noindent
{\bf Constraints.} It is very convenient to enforce simple constraints 
on the state estimates $x_k$. If we take 
\(
\rho_s(x) = \delta_X(x) 
\)
then the prox operator $\prox_{\rho_s}$ is simply the projection onto the 
set $X$. Box constraints are a very common type of constraints 
that enforce known bounds on the state, and have a trivial projection. The proposed framework allows us to use any convex region that has a computationally efficient projection.

\begin{table}
\blue{
\caption{\label{table:PLQ} Prox operators of common PLQ penalties.}
\begin{tabular}{|c|c|c|}\hline
{\bf Penalty $f$} & {\bf $\prox_{\alpha f}(z)$} & {\bf Ref.}\\ \hline
$\frac{1}{2}\|x\|^2$, Fig.~\ref{fig:quadratic} 
& $ \frac{1}{1+\alpha}z$ 
& \cite{freedman2009statistical,SeberWild2003}\\
\hline
$\|x\|_1$, Fig.~\ref{fig:1norm} 
&$ \mbox{sign}(z)\odot(|z|-\alpha)_+$ 
& \cite{Hastie90,LARS2004}\\
\hline
$q_\tau$, Fig.~\ref{fig:quantile} 
&$\begin{cases} z_i-\alpha(1-\tau) & \quad z_i > \alpha(1-\tau)\\
z_i + \alpha\tau & \quad z_i < -\alpha\tau \\
0 & \quad \mbox{else}
\end{cases}$ 
& \cite{KB78,KG01}\\
\hline
$h_\kappa$, Fig.~\ref{fig:huber} 
&$\frac{\alpha}{\alpha + \kappa}z + \frac{\kappa}{\alpha + \kappa} \prox_{(\alpha + \kappa)\|\cdot\|_1}(z)$ 
& \cite{Mar}\\
\hline
$q_{\tau,\kappa}$, Fig.~\ref{fig:huberQ} 
&$\frac{\alpha}{\alpha + \kappa}z + \frac{\kappa}{\alpha + \kappa} \prox_{(\alpha + \kappa)q_\tau}(z)$ 
& \cite{aravkin2014qh}\\
\hline
$v_\epsilon$, Fig.~\ref{fig:vapnik} 
&$\begin{cases} z_i-\alpha & \quad z_i > \epsilon + \alpha\\
\epsilon & \quad \epsilon < z_i \leq \alpha + \epsilon \\
z_i & \quad -\epsilon \leq z_i \leq \epsilon \\
-\epsilon & \quad -\epsilon-\alpha  < z_i \leq - \epsilon \\
z_i + \alpha\tau & \quad z_i < -\alpha -\epsilon. 
\end{cases}$ 
& \cite{Vapnik98}\\
\hline
hubnik-$\kappa$, Fig.~\ref{fig:sel}
&$\frac{\alpha}{\alpha + \kappa}z + \frac{\kappa}{\alpha + \kappa} \prox_{(\alpha + \kappa)v_\epsilon}(z)$ 
& \cite{chu2001unified,lee2005epsi}\\
\hline
e-net, Fig.~\ref{fig:enet} 
&$\prox_{\frac{\alpha}{1 + 2\alpha} \|\cdot\|_1}\left(\frac{1}{1 + 2\alpha}z\right)$ 
& \cite{EN_2005}\\
\hline
\end{tabular}
}
\end{table}

\section{Navigation Model}
\label{sec:navmodels}

We use a
constant-velocity kinematic model that is 
appropriate for many underwater vehicle applications, where
accelerations are heavily damped and trajectories are often long straight lines (e.g. for transit or survey work). When the attitude is known or changing slowly, the model can be linearized effectively. 
For a vehicle that is well-instrumented in attitude, the uncertainty in position
(and the x-y states in particular) is typically orders of magnitude larger than
the uncertainty in attitude. Thus, in practice, we often simplify the full
nonlinear vehicle process model to track only position states 
$( x, y, z)$, while assuming that
the attitude states $(r, p,h)$ are directly available from the most recent sensor
measurements. To make the model linear, the position and its derivatives are
referenced to the local-level frame.

An effective model must counteract biases, outliers, and data discretization in the IMU data. We develop this model using the elements of the proposed framework. 

\noindent
\textbf{Process model.} To incorporate linear acceleration measurements from
an IMU, we track linear velocities 
and linear acceleration in the state vector:
\begin{equation}
  x_s = [x,y,z,\dot{x},\dot{y},\dot{z},\ddot{x},\ddot{y},\ddot{z}]^\top.
\end{equation} \label{lin_pm}
The linear kinematic process model is given by 
\begin{align}
\dot{x}_s &= 
\underbrace{\left[\begin{array}{ccc}
0 & I & 0 \\
0 & 0 & I \\
0 & 0 & 0 \end{array}\right]}_{\mbox{$F_s$}} x_s + 
\underbrace{\left[\begin{array}{c}
0 \\
I \\
0 \end{array}\right]}_{\mbox{$G_s$}} w_s \label{x_s-dot},
\end{align}
where $w_s \sim \mathcal{N}(0,Q_s)$ is zero-mean
Gaussian  noise. 
The linear process model \eqref{x_s-dot} is  
discretized using a Taylor series: 
\begin{align}
x_{s_{k+1}} &= F_{s_k} x_{s_k} + w_{s_k} \label{disc_s}\\
F_{s_k} &= e^{F_s T} \approx
 \approx  \left[\begin{array}{ccc}
I & IT & \frac{1}{2}IT^2\\
0 & I & IT\\
0 & 0 & I \end{array} \right], \nonumber
\end{align}
where $I$ in~(\ref{disc_s}) denotes the $3\times 3$ identity matrix,
and the higher order terms are identically zero because of the
structure of $F_s$. 
We model the process covariance as if the error were the next term in the Taylor series approximation, 
a technique suggested by~\cite{YAA}. More precisely we set covariance to be the outer product, $\Gamma^T \Gamma$ where
\[
\Gamma = \begin{bmatrix}\frac{1}{3!}IT^3 & \frac{1}{2!}IT^2 & IT\end{bmatrix}
\]
This leads to a rank 3 covariance for a $9\times 9$ matrix for a model that comprises
position, velocity, and acceleration.


Given this covariance structure, the process model will penalize changes in acceleration. As the vehicle travels in a relatively straight line with small corrections, we expect to see acceleration mostly constant with a few small jumps. We use the $\ell_1$ norm for innovations, as it encourages exact fits while simultaneously allowing occasional sudden changes.

\noindent
\textbf{Measurement model.} The inertial measurement unit (IMU) measures linear and angular accelerations 
relative to the physical frame of the
vehicle on which it is mounted, while 
the state tracks linear acceleration  relative to the navigation frame. We obtain the coordinate transformation between these frames using heading, pitch, and roll of the vehicle:
\begin{equation}
  R(\varphi) = R^\top_h R^\top_p R^\top_r,
\end{equation}
where $R_h$, $R_p$, and $R_r$ are given by 
\begin{equation}
\left[\begin{array}{ccc}
      c{h} & s{h} & 0 \\
      -s{h} & c{h} & 0 \\
      0 & 0 & 1 \end{array}\right], \quad  
      \left[\begin{array}{ccc}
      c{p} & 0 & -s{p} \\
      0 & 1 & 0 \\
      sp & 0 & cp \end{array}\right], \quad
      \left[\begin{array}{ccc}
      1 & 0 & 0 \\
      0 &c{r} & s{r} \\
      0 & -s{r} & c{r} \end{array}\right]
\end{equation}
with $c\cdot$ and $s\cdot$ shorthand for $\cos(\cdot)$ 
and $\sin(\cdot)$. 

Position data from the USBL is sampled at a lower update rate than the IMU. For any $s_k$ where such position data is available, we have the measurement model 
\[
H_{s_k} = \begin{bmatrix} I_{3 \times 3} & 0_{3 \times 6}\\ 0_{3 \times 6} & R(\varphi)\end{bmatrix}, \quad
z_{s_k} = \begin{bmatrix} \mathrm{usbl}^\top & \ddot{x}_{\textrm{meas}} &  \ddot{y}_{\textrm{meas}} & \ddot{z}_{\textrm{meas}}\end{bmatrix}^\top
\]
If there is no position data measured at time $s$ then we use the model 
\[
H_{s_k} = \begin{bmatrix} 0_{3 \times 3} & 0_{3 \times 6}\\ 0_{3 \times 6} & R(\varphi)\end{bmatrix}, \quad 
z_{s_k} = \begin{bmatrix} 0& \ddot{x}_{\textrm{meas}} & \ddot{y}_{\textrm{meas}} & \ddot{z}_{\textrm{meas}}\end{bmatrix}^\top.
\]
The covariance used for measurement data similarly depends on whether there is position data available:
\[
R_{s_k} = \begin{bmatrix}  U_s & 0_{3 \times 3}\\ 0_{3 \times 3} & r_s I_{3 \times 3}\end{bmatrix}, \quad R_{s_k} = \begin{bmatrix} 0_{3 \times 3} & 0_{3 \times 3}\\ 0_{3 \times 3} & r_s I_{3 \times 3}\end{bmatrix}
\]
where $U$ is a diagonal matrix reflecting position uncertainty, while $r_s$ captures uncertainty in IMU measurements.

\noindent
\textbf{Bias model.}
To compensate for the bias in acceleration data we augment the state vector to include bias variables:
\[
  \bar{x}_s = [x_s^T, b_1, b_2, b_3]^\top
\]
where $b_1,b_2,b_3$ are bias terms for acceleration in the $x,y,z$ directions 
in the local frame. 
To pass the bias estimates forward in time the process matrix is augmented with an identity block.
\begin{align}
    \bar{x}_{s_{k+1}} &= \bar{F}_{s_k}x_{s_k} + \bar{w}_{s_k}\\
    \bar{F}_{s_k} &= \begin{bmatrix}F_{s_k} & 0\\ 0 & I\end{bmatrix}
\end{align}
At the first time point ($s_k = 1$) we augment the covariance matrix with an identity block, 
and at all other time points we augment with a zero block. This adds equality constraints for the bias terms over all time points.
The approach generalized easily to model piecewise-constant biases over longer periods. 

\[
\bar{Q}_{1} = \begin{bmatrix} \Gamma^T\Gamma & 0\\ 0 & I\end{bmatrix}, 
\quad \bar{Q}_{s_k} = \begin{bmatrix} \Gamma^T\Gamma & 0\\ 0 & 0\end{bmatrix} \quad s_k > 1.
\]

The measurement matrices are augmented with an identity block that shifts the acceleration measurements using the bias.

\noindent
\textbf{Discretization model.}
The measurement loss function is chosen to account for the level of discretization present in the data. The Vapnik loss (Figure~\ref{fig:PLQ}f) has a `deadzone' around the origin where small discrepancies are not penalized. This region is set using the level of discretization in the data, in this case $0.05$. The sharp corners of the Vapnik loss encourage errors to lie exactly in them, an unnecessary artifact. Thus we use the `Huberized' version dubbed `hubnik' (Figure~\ref{fig:PLQ} g). 
The prox operators for both losses are computed in Table~\ref{table:PLQ}.

\section{Analysis of Mooring Data}
\label{sec:analysis}
We are interested in the ability to maintain an accurate position
estimate on-board an autonomous underwater vehicle in real-time using acceleration
measurements from a low-cost IMU, given periodic
position fixes.

 To test this, we use the singular robust Kalman framework to analyze data collected from a surface mooring
equipped with an IMU on the subsea node. 
The mooring node, which is drifting with the current, is used as a proxy for a slowly
moving underwater vehicle subject to unknown disturbances. 
We look at the position uncertainty and error accrued over time between the
periodic, world-referenced position fixes provided by the USBL system.

\noindent
\textbf{Data description.}
Position fixes are available from a ship-based Sonardyne Ranger 2 USBL every 2 seconds, which we subsample to varying degrees for the analysis. Linear acceleration data from an LSM303D 3-axis accelerometer was collected at $\sim0.075$ m/s$^2$ precision using a Raspberry Pi Zero. 

\begin{figure}[h!]
\begin{center}
\includegraphics[scale=.45]{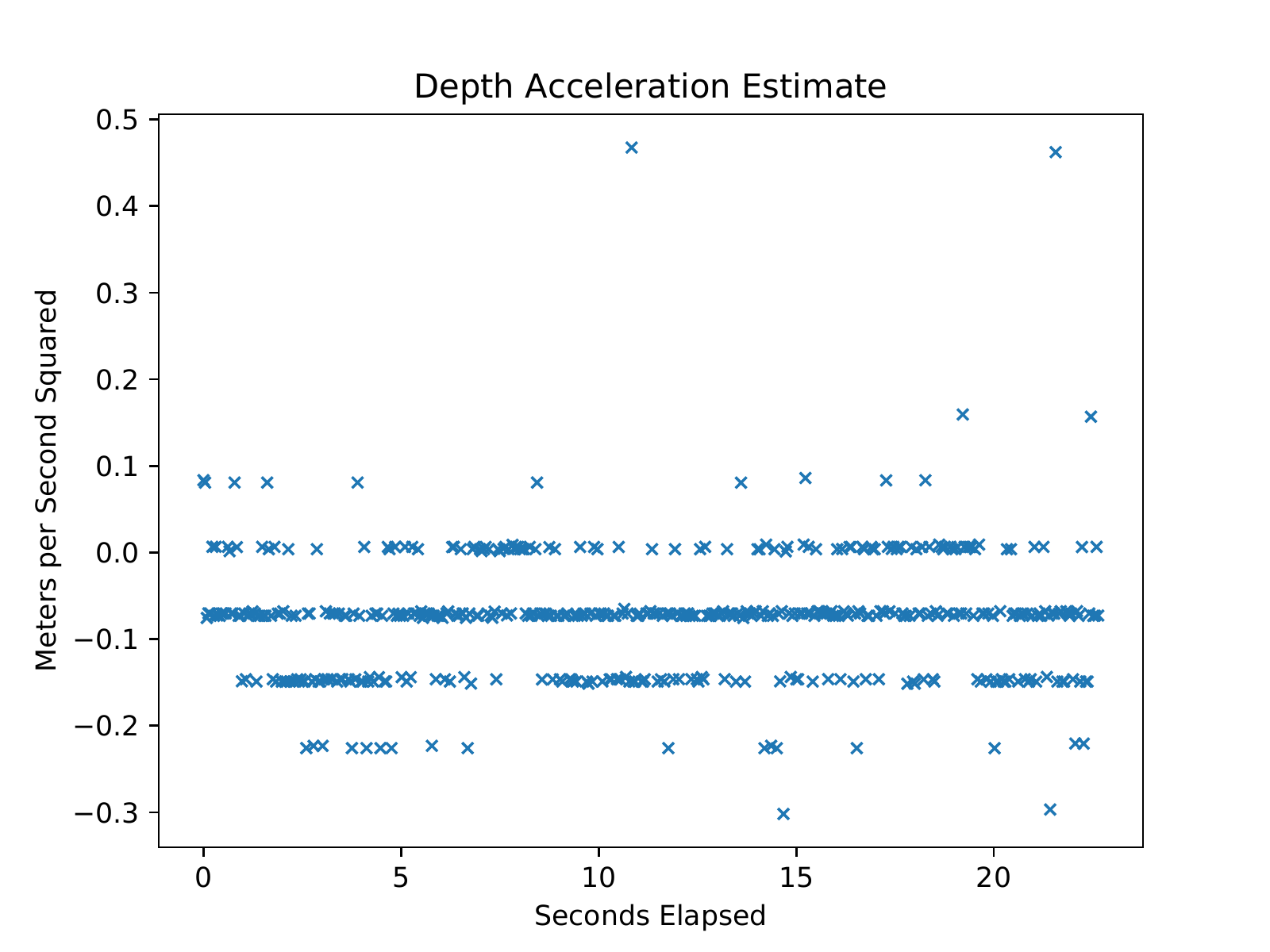}
\caption{A snippet of the depth acceleration data, rotated into the world frame,
  shows the discretization and bias of the acceleration data.}
\label{fig:acceldata}
\end{center}
\end{figure}
\begin{figure}[h!]
\begin{center}
\includegraphics[scale=.45]{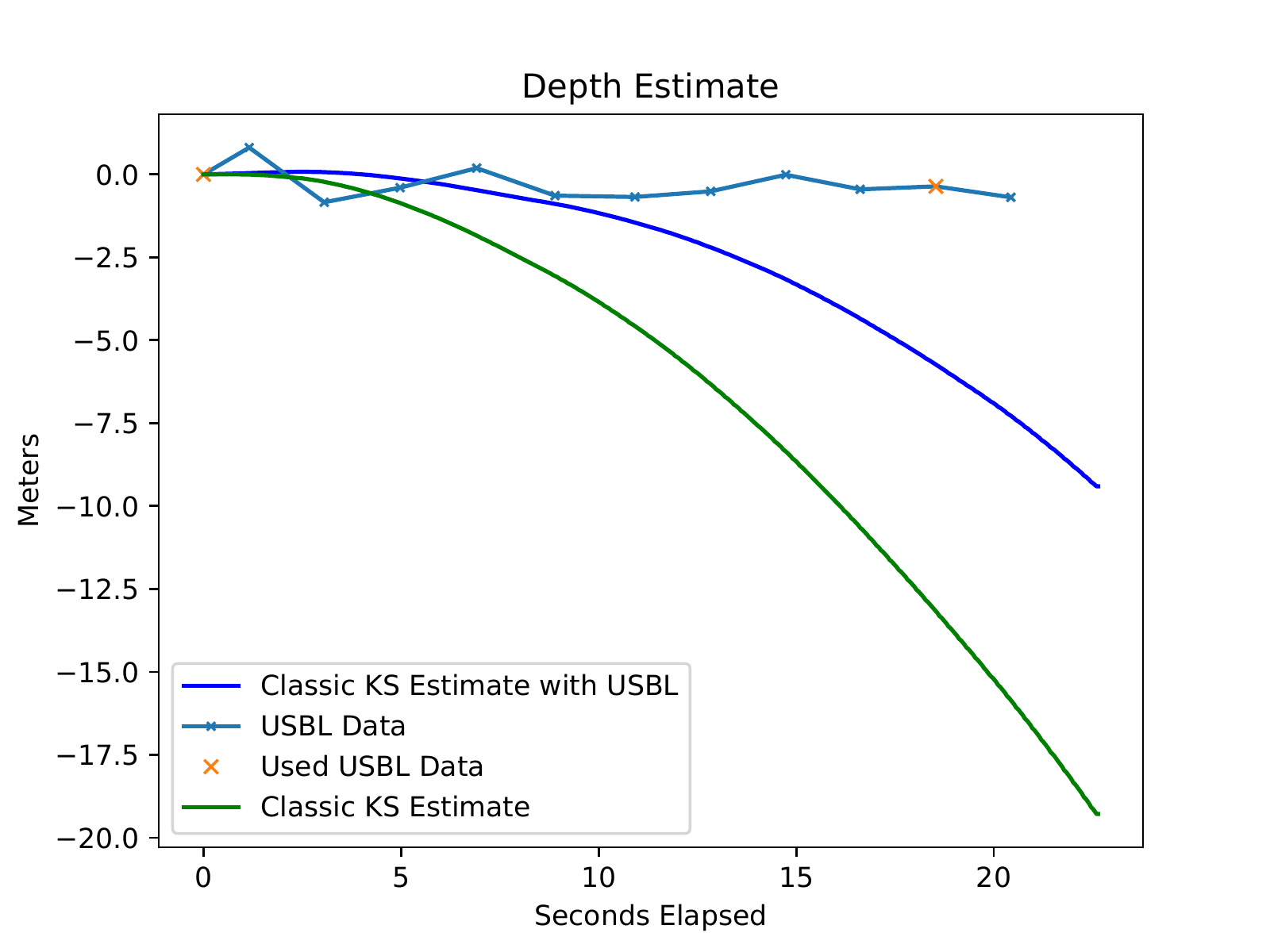}
\includegraphics[scale=.45]{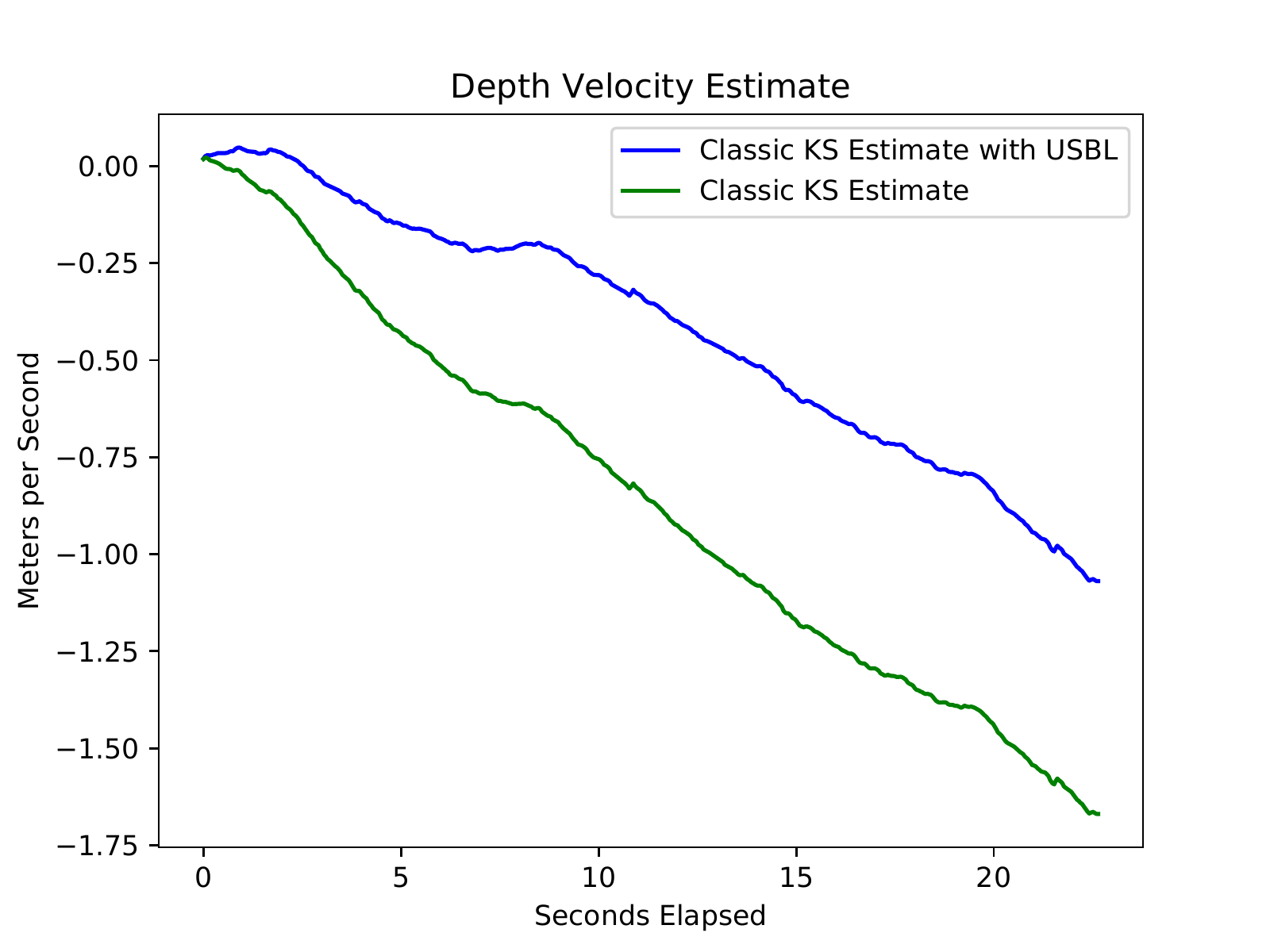}
\caption{Comparison of classic Kalman smoothing applied to depth data with and without USBL position fixes.}
\label{fig:l2_depth}
\end{center}
\end{figure}

\begin{figure}[h!]
\begin{center}
\includegraphics[scale=.45]{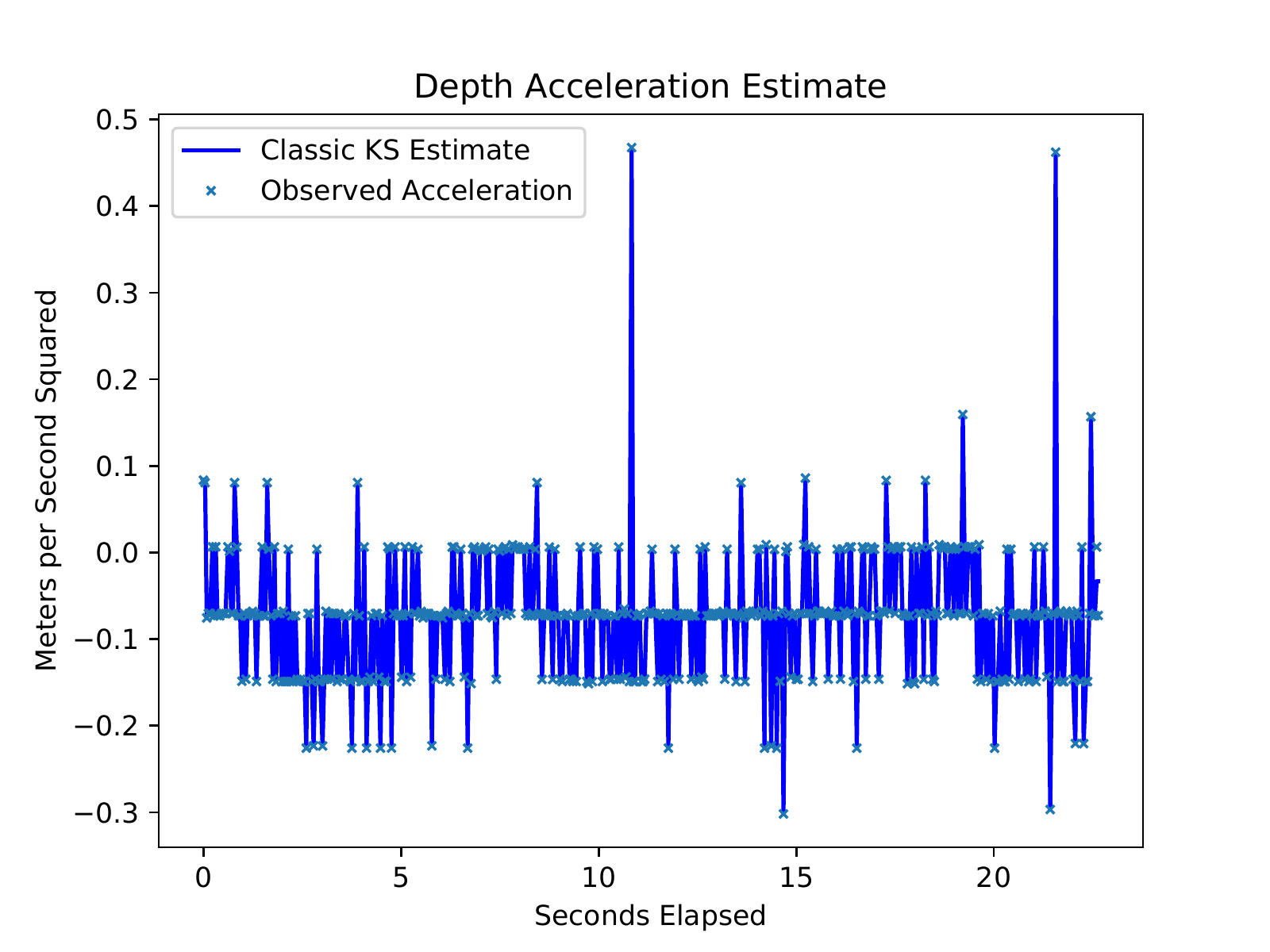}
\includegraphics[scale=.45]{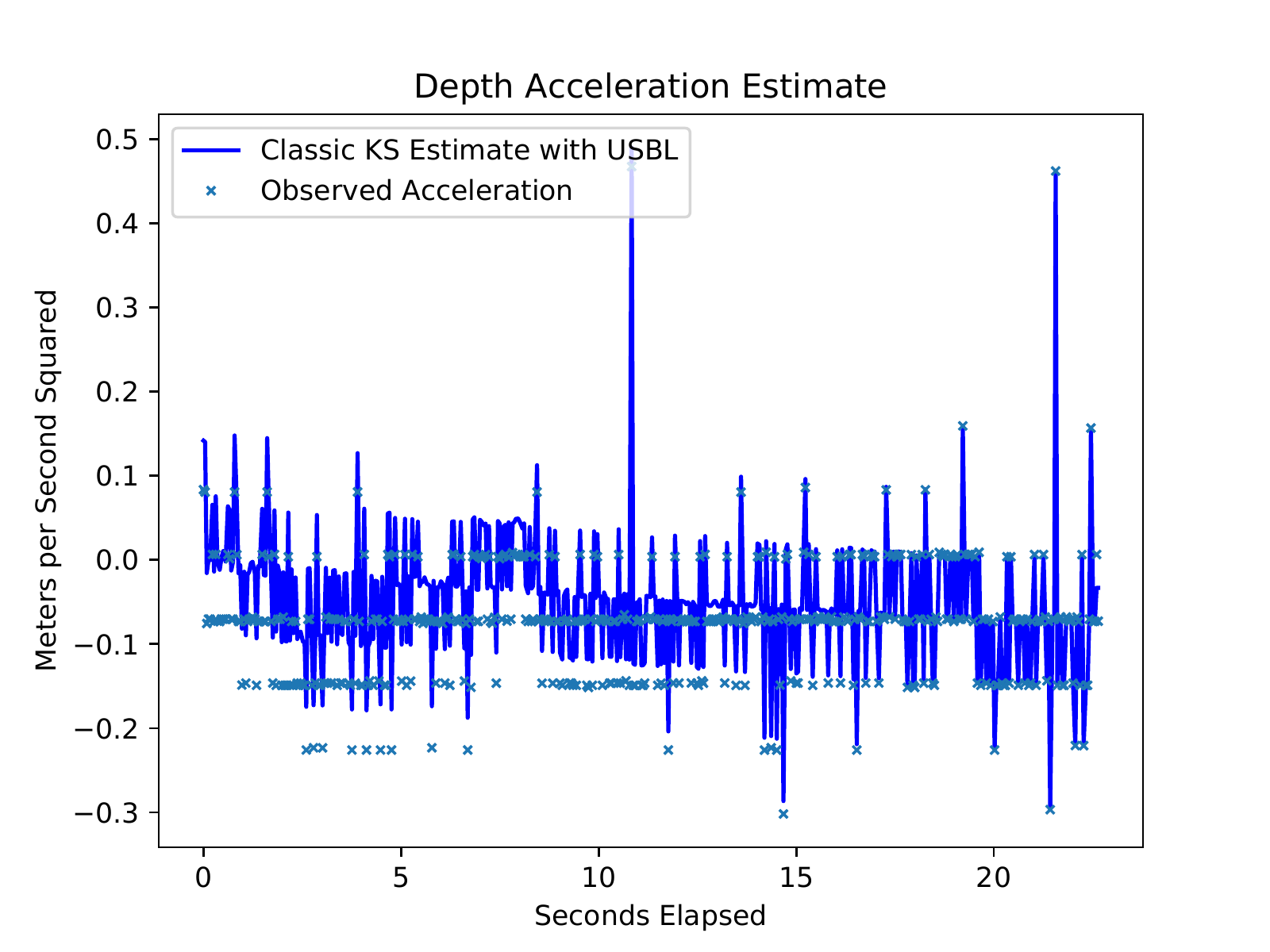}
\caption{Acceleration estimates of classic Kalman smoother without USBL data (top) and with USBL data (bottom).}
\label{fig:l2_accel_depth}
\end{center}
\end{figure}

A small snippet of vertical acceleration data, Figure~\ref{fig:acceldata}, shows the relatively coarse discretization of the measurements along with a mean that is shifted away from zero indicating bias. Outliers also appear likely. Because the IMU was generally upright with small perturbations in pitch and roll, the discretizations in the z-axis of the instrument frame are still clearly visible when the data are rotated into the vertical world frame. The small perturbations in pitch and roll cause the slight variations visible within each discretization level.

To illustrate our model's usefulness on a data set such as in Figure~\ref{fig:acceldata}, we isolate a 25-second window of the data and add modeling elements one at a time, noting the improvements they provide. We finish by running the full model with varying amounts of infrequent position data on a 10-minute section of data as a more practical experiment.

\noindent
\textbf{Model elements applied to data.} 
We consider 25 seconds of IMU data and apply a classic Kalman filter (using least squares measurement and process loss) with the singular navigation model detailed in section~\ref{sec:navmodels}. In practice, one would expect an underwater vehicle to be well instrumented in depth, but for illustration purposes we focus on depth and the vertical acceleration measurements. These measurements are the most biased and therefore improvements made by adding modeling elements are most clearly shown.

The 25 seconds of IMU data starts at a USBL measurement. Position is initialized to this starting USBL measurement, but for ease of comparison that position is treated as (0,0,0) and all subsequent USBL measurements are treated as relative offsets. DRS works better when its reasonably initialized. The initial state vector over the smoothing window is populated by setting acceleration to $0$, and then propagating forward the most recent position fix with a significantly damped measurement of the most recent velocity to prevent divergence of the initial vector.



We start by adding bias estimation (for the acceleration measurements). This requires a second USBL position fix, 19 seconds into the 25-second data series. Figure~\ref{fig:l2_depth} shows position and velocity estimates from a classic Kalman smoother applied to depth acceleration data with and without additional USBL position data. As expected due to the quality of the data we see poor performance in both cases although the small amount of additional position data offers slight improvement.

%
%

Figure~\ref{fig:l2_accel_depth} shows acceleration estimates of the classic Kalman smoother with and without additional USBL data. The two USBL measurements affect 
acceleration estimates locally but are  generally overpowered by the vastly more 
complete set of (biased) acceleration data.  
\begin{figure}[h!]
\begin{center}
\includegraphics[scale=.45]{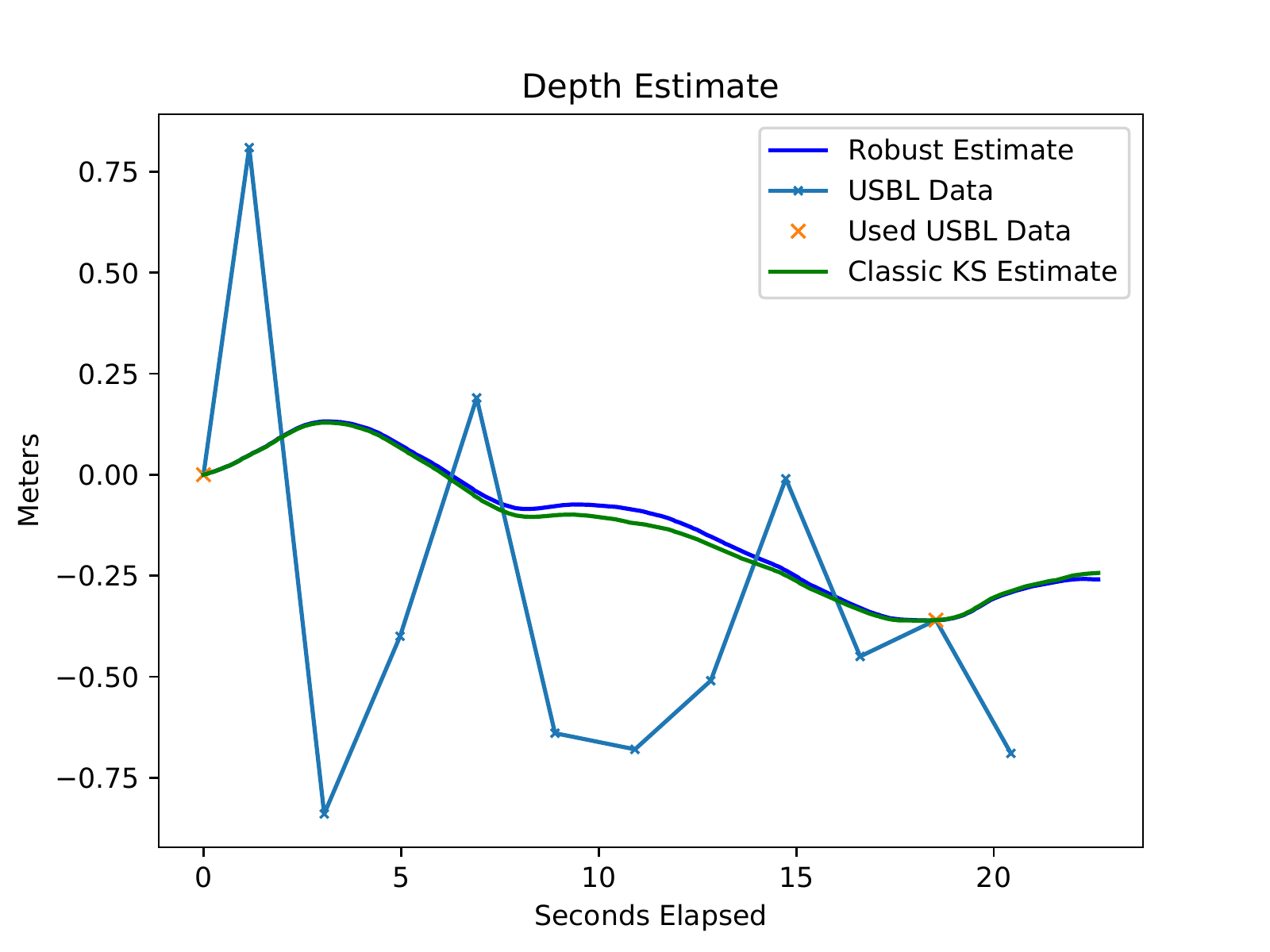}
\includegraphics[scale=.45]{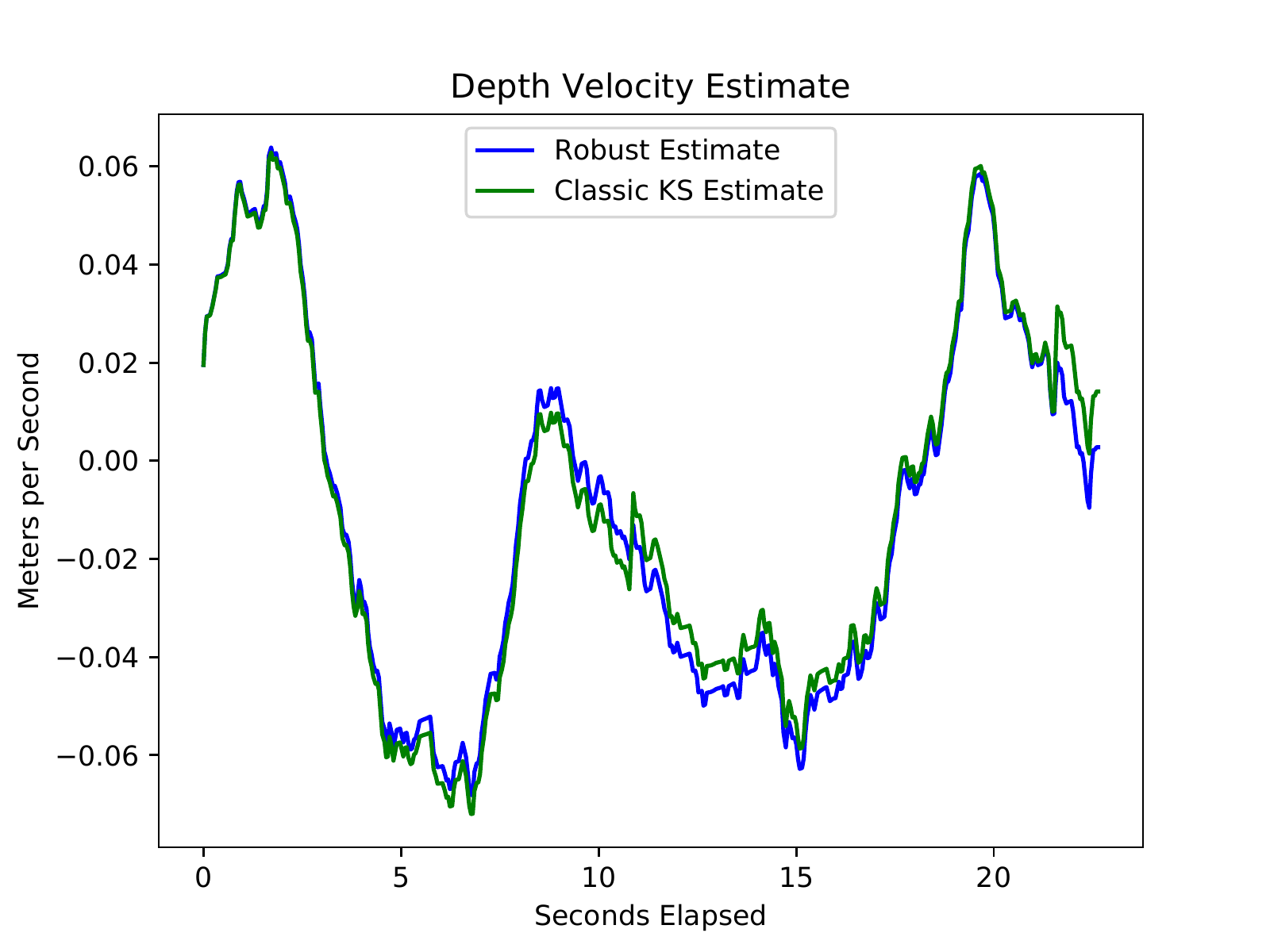}
\caption{Comparison (position and velocity estimates) of debiasing least squares Kalman smoother vs. a robust debiasing smoother equipped with the hubnik loss (Figure~\ref{fig:PLQ}g).}
\label{fig:l2_vs_robust}
\end{center}
\end{figure}

\begin{figure}[h!]
\begin{center}
\includegraphics[scale=.45]{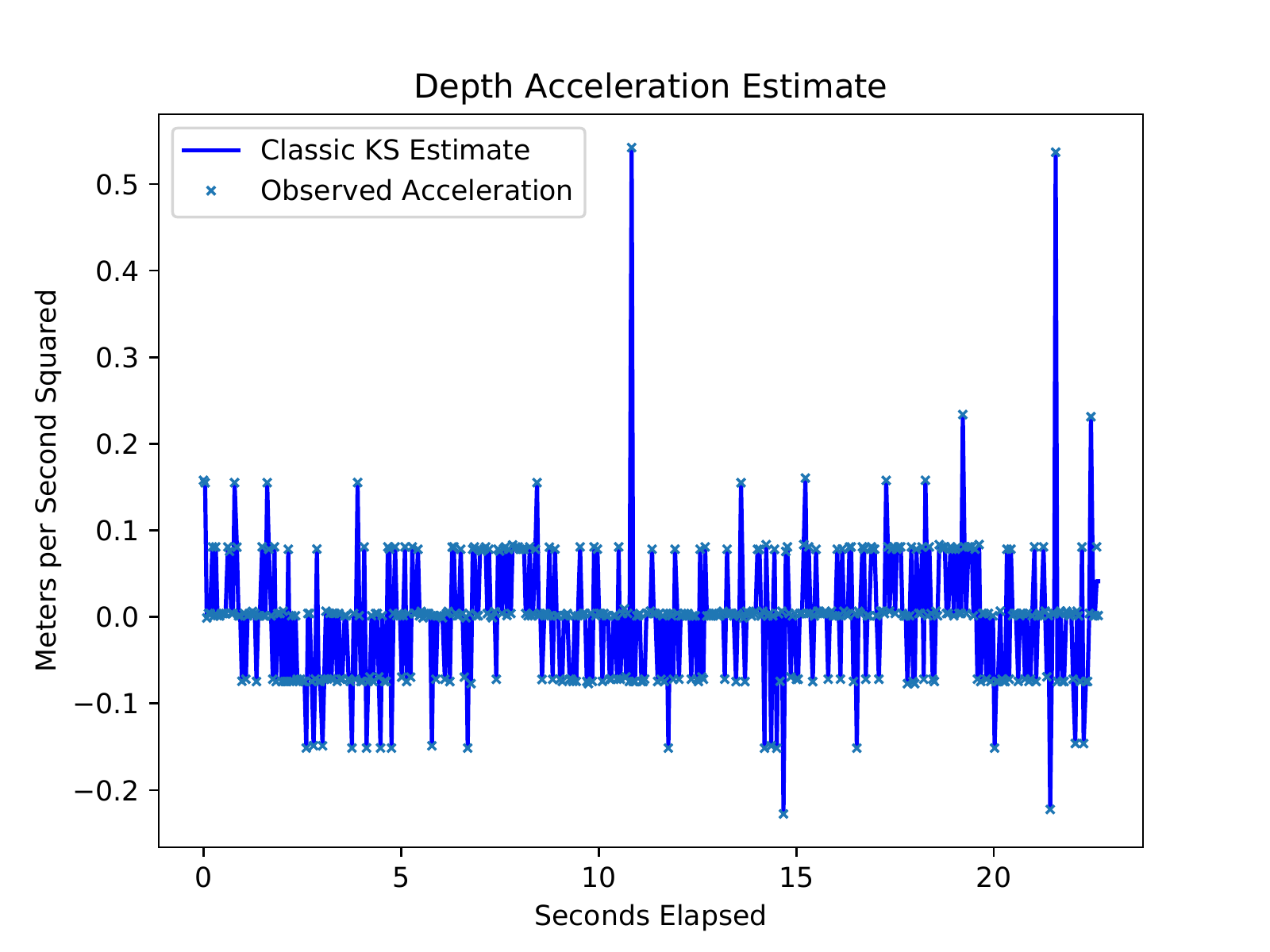}
\includegraphics[scale=.45]{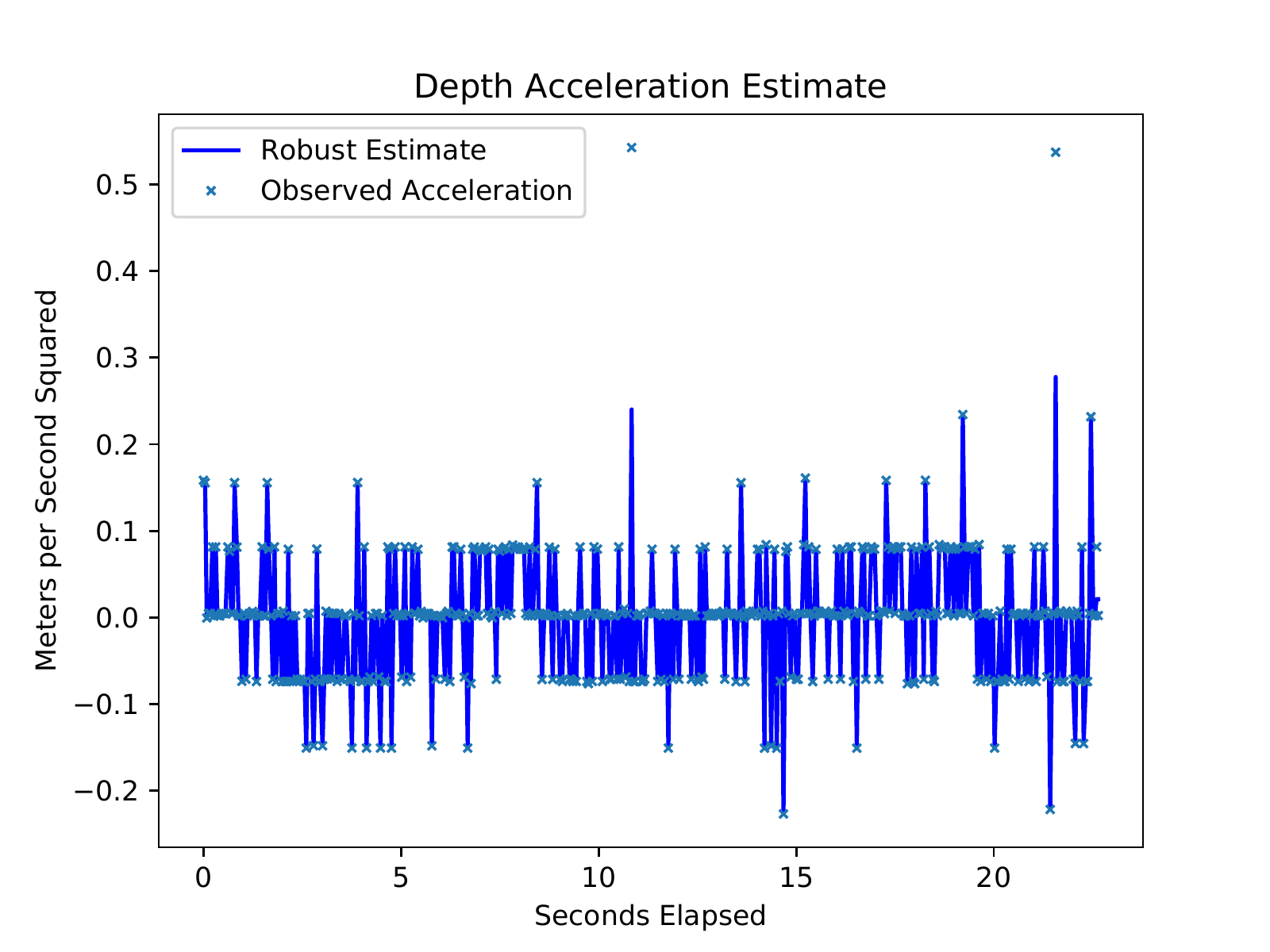}
\caption{Acceleration estimates for least squares debiasing Kalman smoother (top) vs. robust debiasing smoother equipped with the smoothed Vapnik loss (Figure~\ref{fig:PLQ}g) (bottom).}
\label{fig:l2_vs_robust_accel}
\end{center}
\end{figure}

%
The picture improves significantly when we add in bias estimation (via singular process models), using the same two observations (compare Figures~\ref{fig:l2_depth} and ~\ref{fig:l2_vs_robust}). With the bias removed, we can now focus removing 
the noticeable outliers form the acceleration data. For this purpose, we 
use the robust hubnik  (Figure~\ref{fig:PLQ}g) as our measurement loss. The `deadzone' is designed to work with the coarse discretization of the acceleration data.

Figure~\ref{fig:l2_vs_robust} compares the results of (debiased) fitting between the least squares 
Kalman smoother and the robust version using the hubnik loss.
%
%
%
Figure~\ref{fig:l2_vs_robust_accel} shows the acceleration estimates with bias removed in both the classic and robust setting.
There are noticeable differences in the acceleration estimates at around 11 and 22 seconds.
When using the robust smoother, 
the effect of the outliers on the model's estimate is greatly reduced (see
 the velocity estimates at 11 and 22 seconds). 
 In the context of real-time tracking or forecasting, these sudden jumps will yield inaccurate predictions.
 
\noindent
\textbf{10 minutes of data.}
We apply the robust smoother with with bias, outliers, and discretization modeling elements 
to 10 minutes of IMU data. USBL data is available approximately every $2$ seconds, but we test performance of the smoothing algorithm at larger gaps, 
with USBL data supplied at $30, 60,$ and $120$ seconds. 

\begin{figure}[h!]
\begin{center}
\includegraphics[scale=.45]{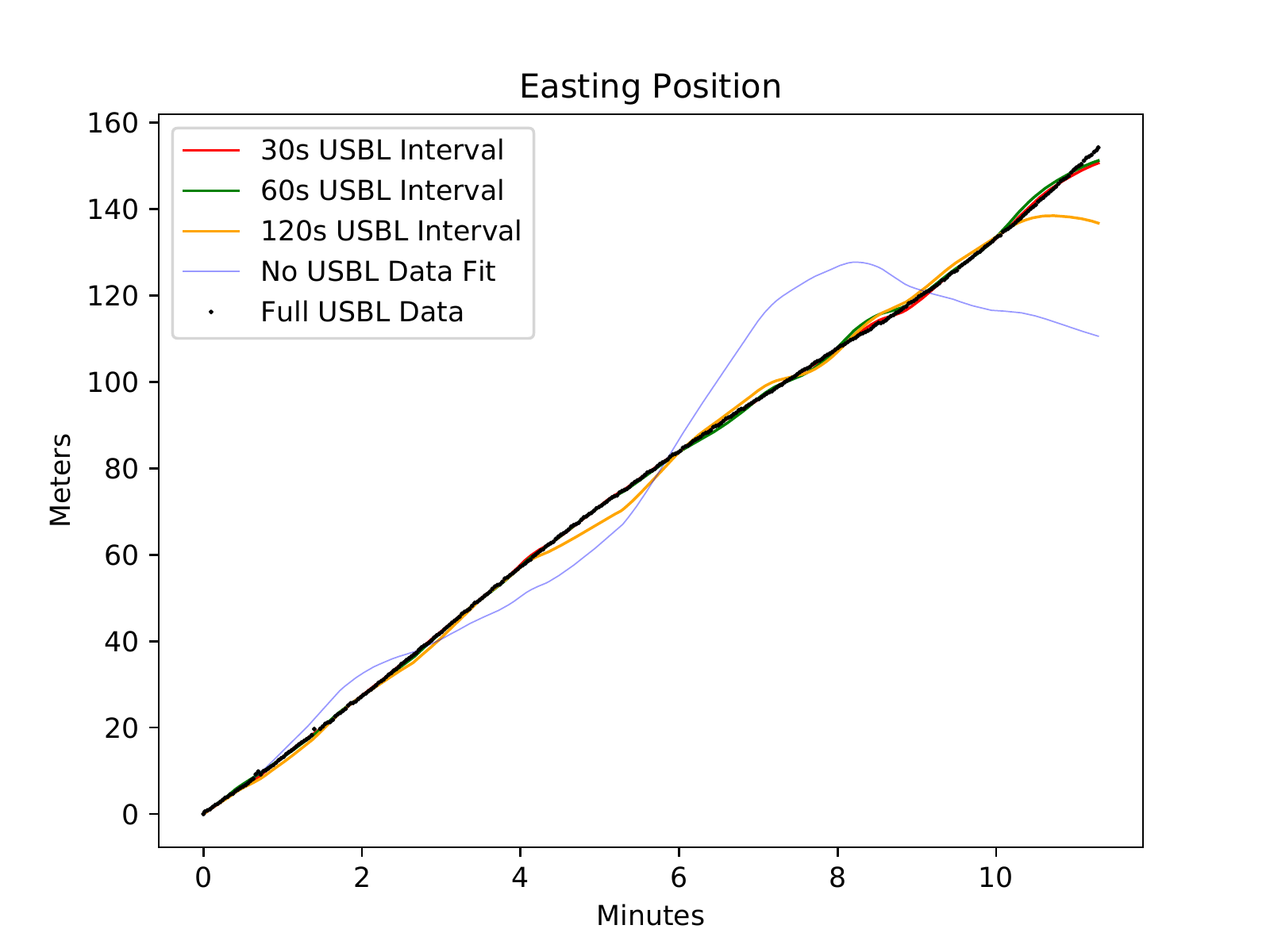}
\includegraphics[scale=.45]{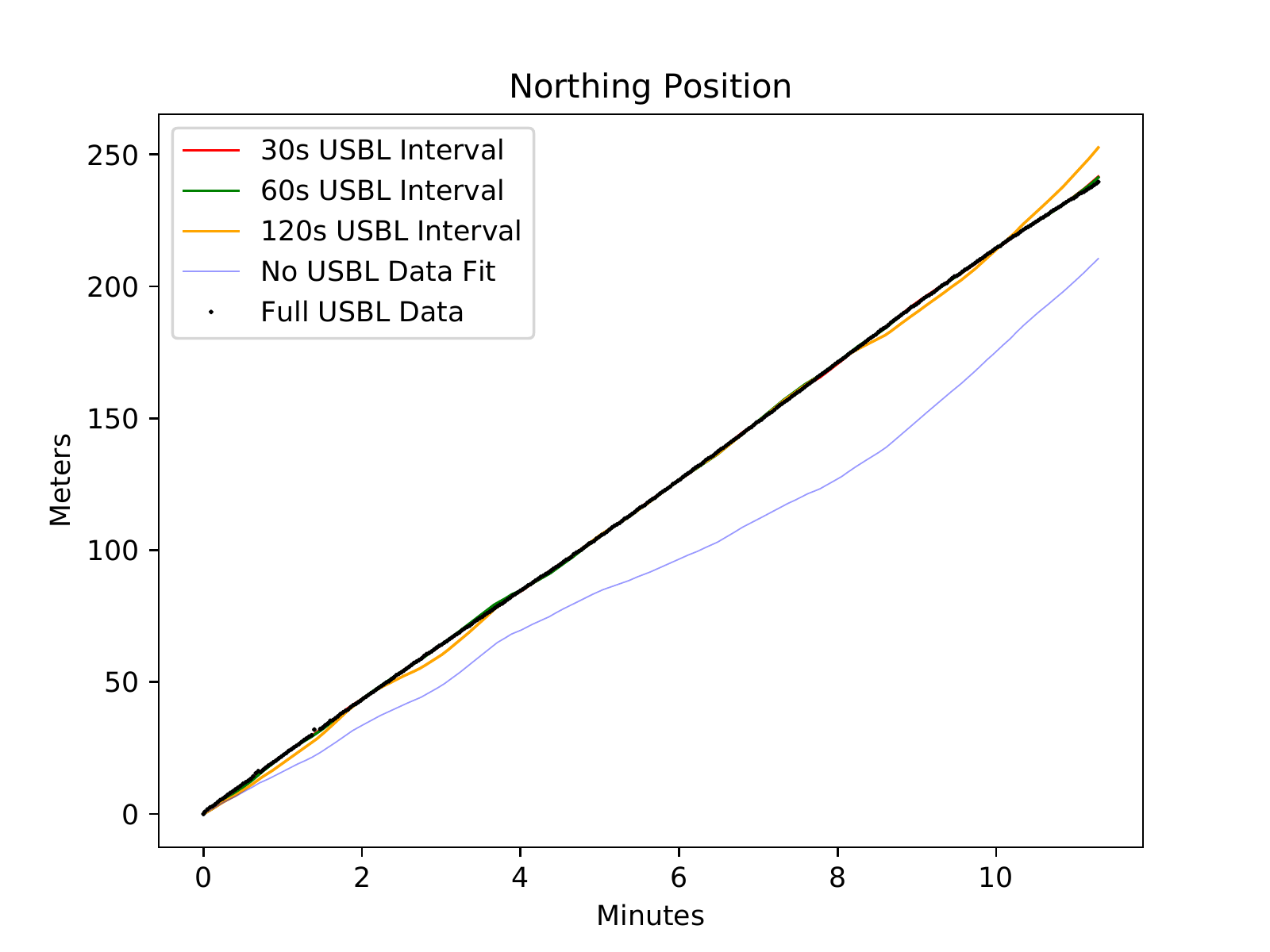}
\caption{Position estimates obtained with robust debiasing smoother for three frequencies of USBL fixes. Robust 
smoothing allows reasonable tracking from infrequent USBL observations.}
\label{fig:posshort}
\end{center}
\end{figure}
\begin{figure}[h!]
\begin{center}
\includegraphics[scale=.45]{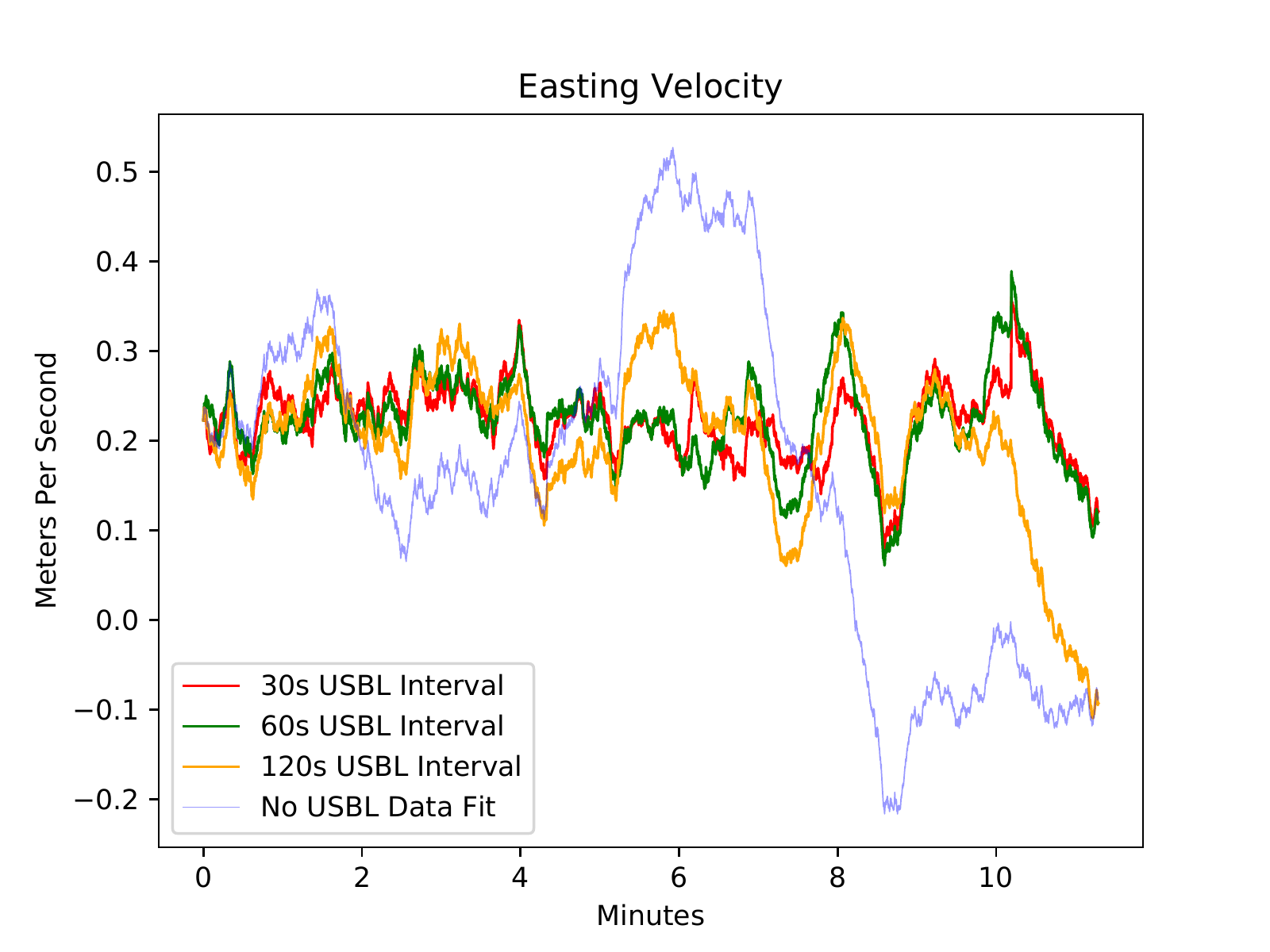}
\includegraphics[scale=.45]{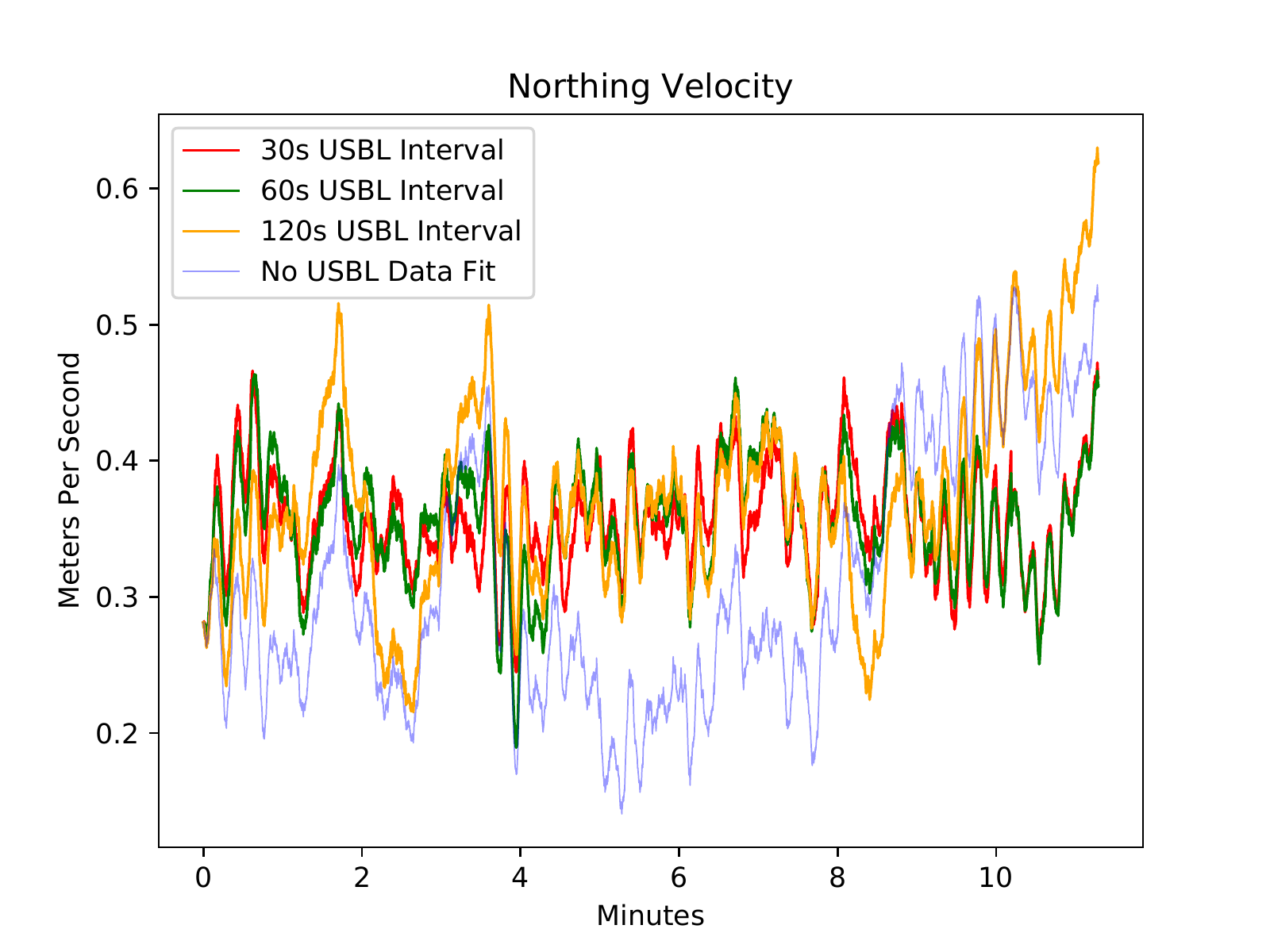}
\caption{Velocity estimates obtained with robust debiasing smoother for different frequencies of position data.  Errors from acceleration measurements build up without USBL fixes, but infrequent USBL measurements still allow velocity estimation.}
\label{fig:velshort}
\end{center}
\end{figure}

Figure~\ref{fig:posshort} has the fitted position plots for all three frequencies. 
We see that without any USBL fix data, the estimate suffers; but even 
infrequent fixes give significant improvements when using the full capability of the smoother. 
In fact there are diminishing returns in increasing the USBL frequency; this is a promising 
result towards the future goal of a practical online implementation, particularly 
in settings where high-quality USBL observations are unavailable (e.g. during dives). 

Figure~\ref{fig:velshort} shows the fitted
velocity. Here the effect of additional USBL data are more apparent, as velocity is completely 
inferred from position and acceleration. However, the smoothing estimates of velocity
at infrequent USBL fixes are still very good compared to those informed by frequent USBL fixes.

\section{Discussion and Future Work}
\label{sec:conclusion}
We have proposed a singular Kalman smoothing framework that can use 
singular covariance models for process and measurements, convex robust 
losses, and state-space constraints. The modeler can use any convex 
loss that has an implementable prox, a class that includes the most common choices 
used for inference in tracking and navigation. 
The framework offers a wide range of flexibility that can be used to either counteract undesirable characteristics present in data or be used to increase model performance based off of relevant field knowledge. Future work will consider real-time implementation, as well as extension to nonlinear models. 

Numerical experiments show that these modeling elements yield significant improvement on a noisy, challenging dataset. We also see that having 
a robust model makes the smoother less reliant on frequent high-quality position updates, 
which is a very promising development for underwater navigation. 


This paper develops several tools required to move to robust singular 
tracking in real-time. A promising aspect of singular noise models 
is that they make it possible to do simple robust windowed smoothing, 
where estimates are constrained between windows as the tracking proceeds. 
Constraints on the state may play a bigger role in real-time estimation, since they can help 
detect outliers faster. Finally, robust penalties that provide better 
estimates may further improve performance of the DRS algorithm, 
by providing an effective initialization 
for each new window. We will focus on these developments in future work. 

\section*{Acknowledgements}

This material is based upon work supported by the Defense Advanced Research Agency (DARPA) and Naval Information Warfare Center (NIWC) Pacific under Contract No. N66001-16-C-4001. It is approved for public release, distribution unlimited. The views, opinions and/or findings expressed are those of the authors and should not be interpreted as representing the official views or policies of the Department of Defense or the U.S. Government.

\bibliographystyle{IEEEtran}
\bibliography{biblio.blb}


\end{document}